# Unit-to-Plant Stability Shaping of Multi-Electrolyzer ReP2H Plants via Interface Design and Dispatch

Miao Zhang, Yiwei Qiu, *Member*, *IEEE*, Xiaoyu Wang, Linlin Wu, Yi Zhou, *Member*, *IEEE*, Shi Chen, *Member*, *IEEE*, Buxiang Zhou, *Member*, *IEEE*, and Kaigui Xie, *Fellow*, *IEEE*

***Abstract*—Alkaline water electrolysis (AWE) units supplied by insulated gate bipolar transistor rectifiers (IGBT-Rs) may experience oscillations caused by coupling between rectifier control and electrolyzer (ELZ) dynamics. Because this risk varies with unit loading and power allocation, production-oriented dispatch may place a multi-ELZ renewable power-to-hydrogen (ReP2H) plant near or exceed its stability boundary. This paper proposes a stability-oriented framework for control design and plant production dispatch. A three-port admittance model links the ac port, dc link, and electrolysis stack. Unit-level dc-port analysis quantifies the effects of loading, temperature, Buck bandwidth, and dc-link capacitance, while plant-level aggregation evaluates how unit commitment and power allocation affect stability. Results show that higher loading reduces stability, whereas larger dc-link capacitance and higher Buck bandwidth improve it. Under the same plant loading, different power allocations result in different plant-level stability margins, with balanced allocation generally providing a larger margin than concentrated allocation. The plant-level model thus distinguishes the stability margins of admissible schedules. Hardware-in-the-loop (HIL) tests validate these trends and the proposed redistribution rule. The resulting operating regions and dispatch rules can be used to screen unit commitment and power allocation decisions in plant production scheduling.**

## I. INTRODUCTION

### A. Background and Motivation

Large-scale renewable power-to-hydrogen (ReP2H) plants typically comprise multiple rectifier-fed electrolyzers (ELZs), connected in parallel at the point of common coupling (PCC) [1], [2]. These plants are often co-built and located near renewable energy bases and therefore operate under weak-grid conditions with limited short-circuit capacity.

As shown in Fig. 1, a typical industry-rated ELZ includes an alkaline water electrolysis (AWE) stack, a Buck converter, a dc link, and an insulated gate bipolar transistor rectifier (IGBT-R) [3]. The IGBT-R regulates the ac current and dc-link voltage, while the Buck converter controls the stack current. Their dynamics are coupled through the dc link, so stack current variations can affect the ac response [4], [5]. A similar cascaded structure appears in railway traction systems, comprising a grid-connected converter, a regulated dc link, and a controlled downstream load, where fast converter dynamics and dc-link coupling can interact with network impedance and cause low-frequency oscillations [6], [7]. Such observations imply comparable oscillatory risks in ReP2H plants.

Beyond converter control, the AWE stack introduces operating-point-dependent dynamics on the dc side [8]. Its voltage response is governed by activation polarization, ohmic loss, and electric-double-layer (EDL) effects, which vary with temperature, current density, and degradation [9], [10]. Variations in stack current pass through the Buck converter and perturb the dc link, while the IGBT-R responds through dc-link voltage and ac current control. This feedback interaction may also cause dc-port oscillations that propagate to the ac side. Since stack dynamics vary with operating conditions, the stability characteristics depend on the operating point and are difficult to assess.

The assessment of stability becomes more complex in large-scale ReP2H systems, where the energy management system (EMS) determines the operating states and power references of individual ELZ units [11]. Different production schedules can place the ELZs at different operating points under the same total production command, changing the plant-level dynamic response and stability margin. Recent small-signal stability studies show that control bandwidth, unit loading, stack temperature, gas pressure, and grid strength can shift system modes and stability boundaries [12], [13]. However, most studies address a single unit or homogeneous aggregation and leave the combined stability of ELZs at different operating points unresolved. A port-based framework is thus needed to retain dc-link coupling, identify unit stability boundaries, and formulate stability constraints for unit commitment and power allocation.

### B. Literature Review and Research Gap

Large-scale ReP2H plants commonly comprise multiple rectifier-fed ELZs, with each rectifier supplying one unit. The plant at the *China Energy Engineering Songyuan Hydrogen Industry Park* [14] employs 32 thyristor rectifiers and 32 IGBT rectifiers to supply 64 AWEs, and the *Wind Solar Hydrogen Ammonia Integration Project* [15] in Da'an uses 12 thyristor rectifiers and 24 IGBT-Rs for 36 AWEs. Zeng *et al*. [16] optimized the portfolio of rectifiers. Gao *et al*. [17] and Koponen *et al*. [18] compared their techno-economic performance, while Meng *et al*. [19] proposed a hybrid rectifier to improve power quality. These works support rectifier selection and steady-state design; however, the dynamic interaction among parallel IGBT-R-fed ELZs has received limited attention [20].

Multi-ELZ operation has been widely studied from the scheduling perspective. Zeng *et al*. [21] proposed coordinated active and reactive power dispatch. Qiu *et al*. [22] included thermal and impurity constraints in plant scheduling. Li *et al*. [23] developed an energy management method for grid-connected

Financial support came from the National Natural Science Foundation of China (52377116 and 52312076). (Corresponding author: Yi Zhou)

M. Zhang, Y. Qiu, Y. Zhou, S. Chen, B. Zhou, and K. Xie are with the College of Electrical Engineering, Sichuan University, Chengdu 610065, China. (yizhou3230@scu.edu.cn)

X. Wang and L. Wu are with the State Grid Jibei Electric Power Co., Ltd. Research Institute, Beijing, China, and the North China Electric Power Research Institute Co., Ltd., Beijing, China. L. Wu is also with the State Grid Jibei Zhangjiakou Wind-Solar-Storage-Transmission New Energy Co., Ltd., Zhangjiakou, China.

multi-stack plants. Firdous *et al*. [24] studied the flexibility of large-scale P2H systems. These methods optimize unit commitment, power allocation, and state transitions under power balance and process constraints [25]. Although these decisions change the operating point and dynamic response of each unit, electrical stability has rarely been considered.

The dc-side dynamics of AWE units have been studied through electrochemical and converter models. Ulleberg [26], Varela *et al*. [27], and Iribarren *et al*. [28] developed AWE models for system analysis and dispatch. Guilbert *et al*. [29] and Guo *et al*. [30] studied dc converters and power electronic interfaces for electrolysis. These studies show that stack parameters, operating conditions, and Buck control affect unit dynamics. However, the reported analyses mainly characterize selected operating conditions and do not determine the boundary between stable and unstable states in grid-connected operation.

Impedance-based methods provide a suitable means of determining this boundary because they quantify stability from interactions among connected ports [31]. Pedra *et al*. [32] and Zong *et al*. [33] developed three-port and multiport models for common voltage-source converters (VSCs) that retain ac and dc coupling. Fang *et al*. [34] assessed weak-grid stability of a single-unit rectifier-interfaced hydrogen system. These existing studies mainly consider generic converters or represent the hydrogen plant as a single unit. They do not link the dc-port stability of each IGBT-AWE unit with the aggregated ac-port response of units operating at different points. Consequently, the plant-level stability region associated with unit commitment and power allocation remains unavailable. Retaining the *d*-axis and *q*-axis ac ports together with the dc-link port is necessary because the dc-port termination changes each unit's ac-port response before aggregation at the PCC. The three-port framework therefore supports both unit-level dc-port boundary derivation and plant-level stability assessment without discarding the ac/dc coupling.

### C. Contributions of This Work

To address these gaps, this paper develops a stability-oriented operating framework for large-scale ReP2H plants. The main contributions are as follows:

1) An operating-point-dependent three-port admittance model is established for the complete IGBT-AWE unit. The model retains the dc-link as an explicit port, enabling unit-level dc-port reduction and plant-level ac-port aggregation within one formulation. Its accuracy is validated by a frequency sweep against a detailed switching model.

2) For each ELZ, the dc-port stability margin is mapped over stack loading, temperature, Buck bandwidth, and dc-link capacitance to construct the stability region for operating and interface design. A multi-input multi-output (MIMO) Nyquist distance is further introduced to evaluate the plant-level stability margin under different unit-commitment and power allocation schemes in a multi-ELZ plant.

3) A lookup-based stability-screening rule is developed for multi-ELZ dispatch. Analytical results and multi-unit hardware-in-the-loop (HIL) tests show that different power allocations can produce different stability margins and dynamic responses under the same total production command. They also confirm that stability-oriented power redistribution can increase the plant-level stability margin and suppress oscillations.

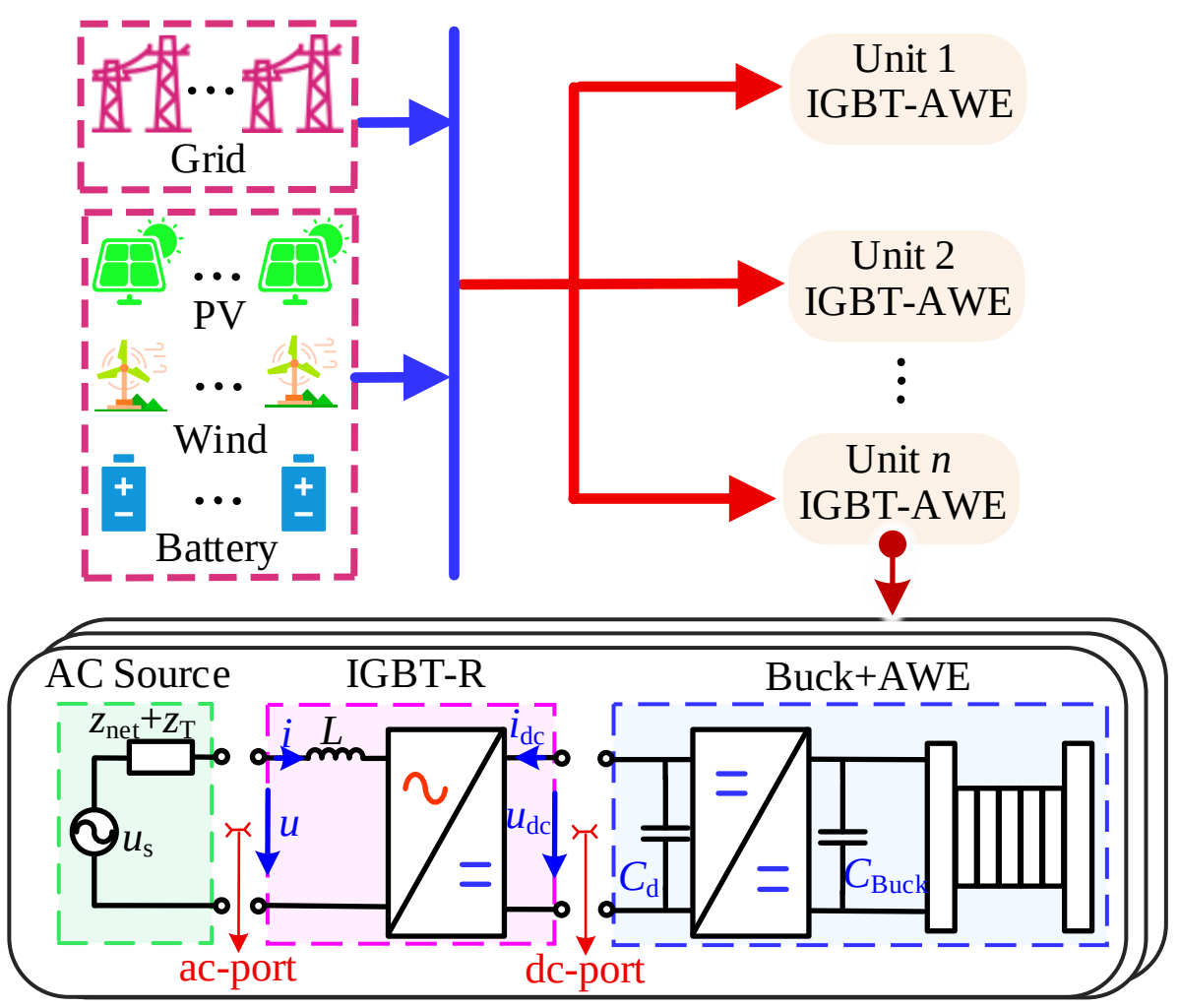

Fig. 1. Topology of a multi-ELZ ReP2H plant and the three-port representation of an IGBT-AWE unit

The rest of this paper is organized as follows: Section II presents the three-port admittance model for the IGBT-AWE units and stability criteria. Section III analyzes the effects of operating conditions and interface parameters on stability. Section IV develops unit-level operating-and-design regions and derives a plant-level stability-oriented power-allocation rule. Section V presents the HIL validation. Section VI concludes the paper.

## II. Port-Based Modeling and Stability Criteria of the Multi-ELZ ReP2H Plant

As shown in Fig. 1, a typical ReP2H plant consists of multiple IGBT-AWE units connected in parallel. Each unit includes an IGBT-R, a capacitive dc link, a current-controlled Buck converter, and an AWE stack. To characterize the coupling between the ac-grid interface and the dc-side electrolysis dynamics, this section establishes a three-port admittance model for a representative unit, with the *dq*-frame ac port and the dc-link port explicitly retained.

### A. Three-Port Admittance Model of the IGBT-AWE Unit

The IGBT-R connects the ac bus to the Buck-interfaced AWE stack through the dc link. A two-port ac representation would eliminate the dc-link termination through which the Buck-AWE dynamics affect the rectifier's ac-port response. The IGBT-R is therefore represented by a three-port admittance model with the *d*-axis, *q*-axis, and dc ports. Linearization around a steady-state operating point gives:

$$\begin{bmatrix}\Delta i_d^s & \Delta i_q^s & \Delta i_{dc}\end{bmatrix}^{\mathrm{T}} = Y_{\text{IGBT-R}}\begin{bmatrix}\Delta u_d^s & \Delta u_q^s & \Delta u_{dc}\end{bmatrix}^{\mathrm{T}}, \quad (1)$$

where

$$Y_{\text{IGBT-R}} = \begin{bmatrix} Y_{dd} & Y_{dq} & Y_{d,dc} \\ Y_{qd} & Y_{qq} & Y_{q,dc} \\ Y_{dc,d} & Y_{dc,q} & Y_{dc,dc}\end{bmatrix}, \quad (2)$$

where the superscript *s* denotes the system *dq*-frame; Δ denotes small-signal perturbation; $\Delta i_d$, $\Delta i_q$ and $\Delta i_{dc}$ are the *d*-axis, *q*-axis, and dc-link current perturbations of the IGBT-R, respectively; $\Delta u_d$, $\Delta u_q$ and $\Delta u_{dc}$ are the corresponding voltage perturbations. In $Y_{IGBT\text{-}R}$, the upper-left 2×2 block represents the ac-port $dq$ admittance. The elements $Y_{d,dc}$ and $Y_{q,dc}$ represent the influence of dc-link voltage perturbation on the ac-port current response, while $Y_{dc,d}$ and $Y_{dc,q}$ characterize the coupling from the ac-port voltages to the dc-port current. The element $Y_{dc,dc}$ is the dc-port self-admittance. The detailed derivation of $Y_{IGBT\text{-}R}$ is provided in Appendix A.

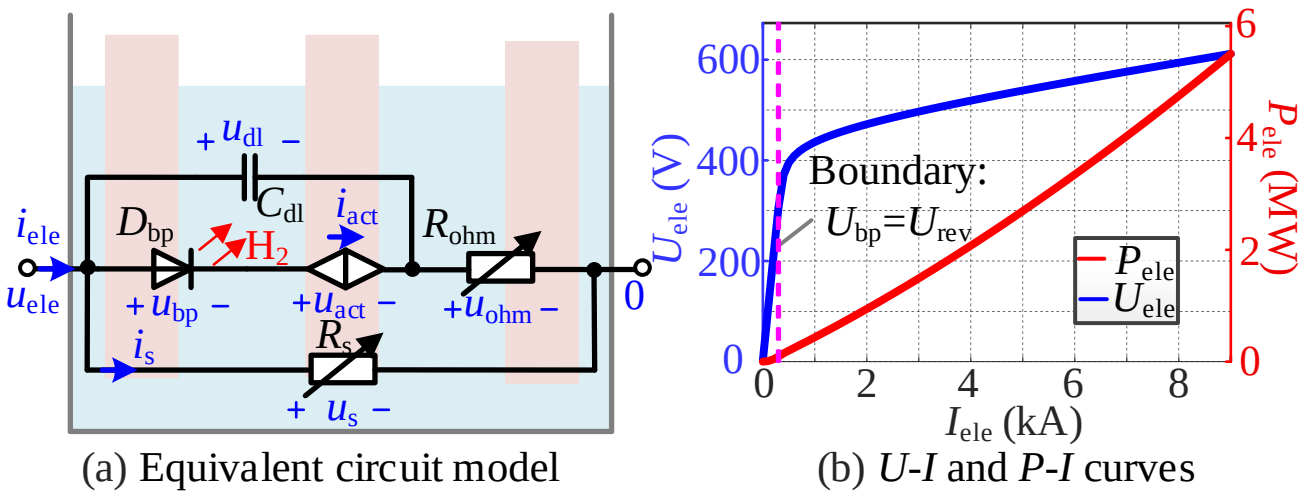

Fig. 2. Equivalent electrical model and steady-state characteristics of the AWE stack

### B. DC-Port Model of the Buck-Interfaced AWE Stack

The dc port of the IGBT-R is terminated by the dc-link capacitor and the Buck-interfaced AWE stack. The corresponding dc-side admittance is expressed as:

$$Y_{dcrc} = sC_d + Y_{Buck\text{-}AWE}, \tag{3}$$

where $C_d$ is the dc-link capacitance and $Y_{Buck\text{-}AWE}$ is the dc-input admittance of the Buck-interfaced AWE stack. The latter is determined by the AWE stack, Buck power stage, and current controller dynamics; the corresponding model is derived below.

As shown in Fig. 2(a), the Buck-interfaced AWE model retains the activation, ohmic, shunt, and EDL dynamics over the frequency range considered in this study. The electrochemical dynamics are formulated at the cell level, with the stack terminal voltage and current defined in the final line of (4), as follows:

$$\begin{cases} u_{cell} = u_{dl} + r_{ohm}(j - j_s),\ C_{dl}\dfrac{du_{dl}}{dt} = j - j_s - j_{act}, \\ j_s = \dfrac{u_{cell}}{r_s},\ j_{act} = f_T(u_{dl}, T),\ f_T = (e^{[\frac{2F}{RT}\alpha(u_{dl}-u_{bp})]} - 1)\beta, \end{cases} \tag{4}$$

where $u_{cell}$ and $u_{dl}$ are the single-cell terminal and EDL voltages, respectively; $j$, $j_s$ and $j_{act}$ are the total, shunt, and activation current densities; and $r_{ohm}$, $r_s$ and $C_{dl}$ are the area-specific ohmic resistance, shunt resistance, and EDL capacitance; the nonlinear function $f_T$ denotes the modified Tafel activation relation [35]; $F$ denotes the Faraday constant; $R$ is the gas constant; $T$ is the temperature; $\alpha$ is the activation coefficient; and $\beta$ is the exchange current density.

For a specified stack-current operating point, (4) is solved to determine the corresponding steady-state stack voltage and internal electrochemical variables. The resulting stack voltage-current characteristic is shown in Fig. 2(b).

For small-signal impedance modeling, the nonlinear activation relation is linearized around a steady-state operating point:

$$\Delta j_{act} = h_1 \Delta u_{dl},\ h_1 = \left.\frac{\partial f_T(u_{dl}, T)}{\partial u_{dl}}\right|_0, \tag{5}$$

TABLE I
PARAMETERS OF A TYPICAL 5-MW-RATED IGBT-AWE UNIT

| Symbols | Parameter | Value |
|---|---|---|
| $u_{d0}$ | $d$-axis steady voltage | 565 V |
| $i_{d0}$ | $d$-axis steady current | 5993 A |
| $u_{dc0}$ | Steady dc-link voltage | 1500 V |
| $f_0$ | Fundamental frequency | 50 Hz |
| $R_f$, $L_f$ | Filter resistor and inductor | 0.0015 Ω, 200 μH |
| $R_g$, $L_g$ | Equivalent source resistor and inductor | 0.0011 Ω, 35 μH |
| $C_d$ | Dc-link capacitor | 18 mF |
| $K_{P_PLL}$, $K_{i_PLL}$ | PLL control | 0.08, 2 |
| $K_{P_ACC}$, $K_{i_ACC}$ | ACC control | 0.075, 2 |
| $K_{P_DVC}$, $K_{i_DVC}$ | DVC control | 15, 20 |
| $L_{Buck}$,$C_{Buck}$ | The inductor and capacitor of the Buck converter | 1 mH, 9 mF |
| $K_{P_Buck}$, $K_{i_Buck}$ | Constant-current control of the Buck converter | $8\times10^{-5}$, $2\times10^{-3}$ |
| $N_{cell}$ | Number of series-connected cells | 242 |
| $A_{cell}$ | Active area of a single cell | 20,000 cm² |
| $U_{rev}$ | Single-cell reversible voltage | 1.228 V |
| $R_{ohm,e}$ | Stack-equivalent ohmic resistance | 15.27 mΩ |
| $R_{s,e}$ | Stack-equivalent shunt resistance | 1.02 Ω |
| $C_{dl,e}$ | Stack-equivalent double-layer capacitance | 15.69 F |
| $F$ | Faraday constant | 96485 C/mol |
| $R$ | Gas constant | 8.31 J/mol·K |
| $\alpha$ | Charge transfer coefficient | 0.1521 |
| $\beta$ | Exchange current density | $1.212\times10^{-4}$ A/cm² |
| $T$ | Stack temperature | 353 K |

where $h_1$ is the area-specific incremental activation conductance. For a stack with $N_{cell}$ series connected cells and an active cell area of $A_{cell}$, the terminal perturbations satisfy:

$$\Delta u_{ele} = N_{cell}\Delta u_{cell},\ \Delta i_{ele} = A_{cell}\Delta j. \tag{6}$$

The corresponding stack-level parameters are then obtained:

$$\begin{cases} R_{ohm,e} = \dfrac{N_{cell} r_{ohm}}{A_{cell}},\quad R_{s,e} = \dfrac{N_{cell} r_s}{A_{cell}}, \\ C_{dl,e} = \dfrac{A_{cell} c_{dl}}{N_{cell}},\quad H_{1,e} = \dfrac{A_{cell} h_1}{N_{cell}}. \end{cases} \tag{7}$$

The stack-terminal small-signal impedance is then directly obtained as

$$Z_{AWE}(s) = R_{s,e} \,||\, \left(R_{ohm,e} + \frac{1}{sC_{dl,e} + H_{1,e}}\right). \tag{8}$$

Equation (8) explicitly retains the operating-point-dependent activation dynamics, the EDL effect, the ohmic voltage drop, and the shunt path, while avoiding separate definitions of intermediate impedances.

The obtained $Z_{AWE}$ is used as the load impedance of the current-controlled Buck converter [36]. Linearizing the Buck power stage and current controller around the same operating point gives:

$$Y_{Buck\text{-}AWE}(s) = \frac{\Delta i_{dc}(s)}{\Delta u_{dc}(s)} = C_b\left(sI - A_b\right)^{-1} B_b + D_b, \tag{9}$$

where $A_b$, $B_b$, $C_b$ and $D_b$ are determined by the Buck inductor, output capacitor, current controller and $Z_{AWE}$. Substitution of

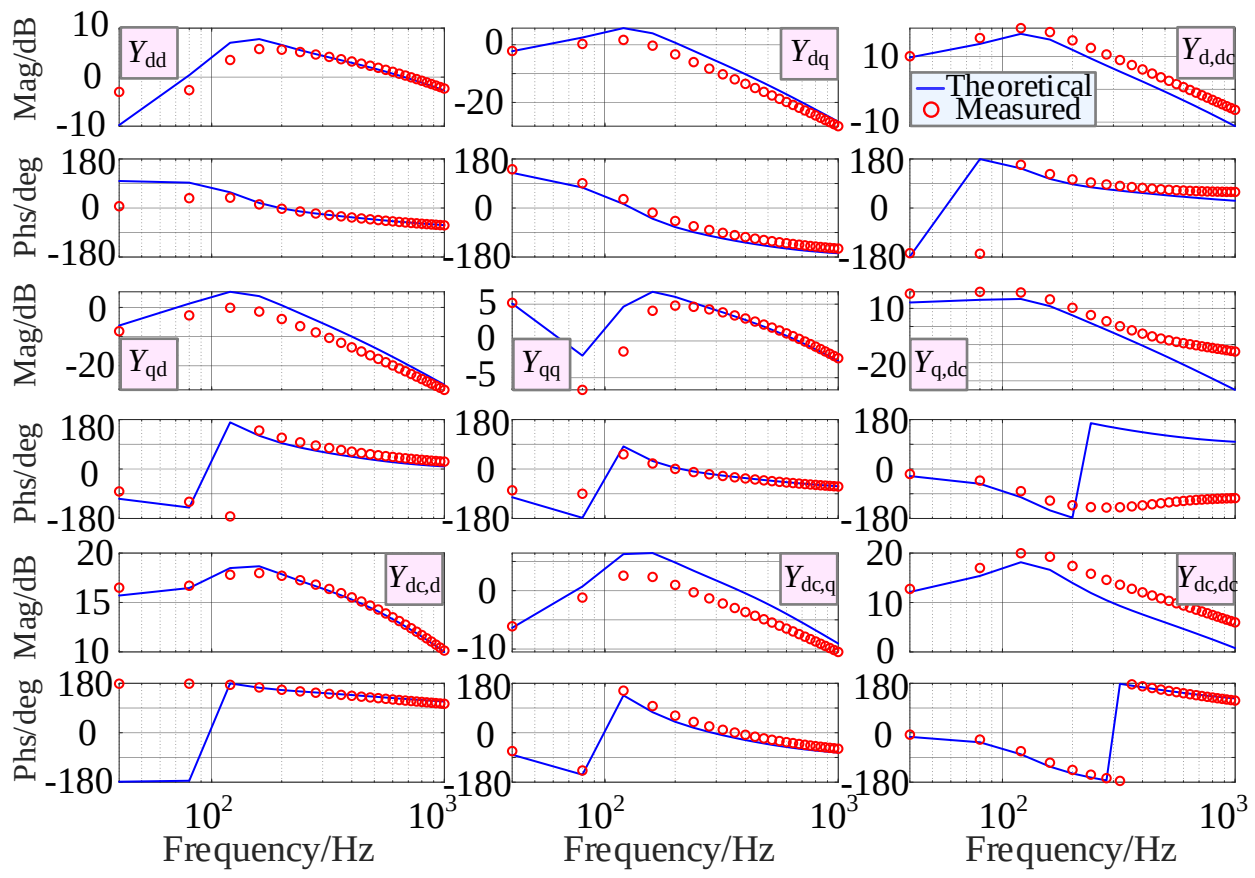

Fig. 3. Frequency-sweep validation of the analytical three-port admittance model against the detailed switching model of the complete IGBT-AWE unit

(9) into (3) gives the dc-side termination admittance $Y_{dcrc}$ used in the subsequent port-reduction and stability analyses.

The parameters of a typical 5-MW-rated IGBT-AWE unit are listed in Table I [37]. The analytical model is validated against a detailed switching model by independent perturbations at the *d*, *q*, and dc ports from 1 to 1000 Hz. Fig. 3 shows that it reproduces the dominant magnitude, phase, and coupling characteristics. The visible deviations occur mainly in weak coupling terms, whose relative errors are amplified by the logarithmic scale and switching ripple. The model is therefore adequate for the subsequent stability analysis.

### C. Unit-Level DC-Port SISO Stability Criterion

The three-port model can be reduced at the dc port to evaluate the internal dc-link interaction of one IGBT-AWE unit [38]. The external network admittance is expressed as:

$$Y_{net} = \begin{bmatrix} Y_{net}^{2\times2} & 0_{2\times1} \\ 0_{1\times2} & Y_{dcrc} \end{bmatrix}, \tag{10}$$

where $Y_{net}^{2\times2}$ denotes the 2×2 ac-network block in the system dq-frame, and $Y_{dcrc}$ is the dc-side admittance defined in (3). The characteristic equation of the IGBT-R connected to the ac network and dc-side admittance is:

$$\det(Y_{net} + Y_{IGBT\text{-}R}) = 0. \tag{11}$$

The matrix in (11) is partitioned according to the ac and dc ports as:

$$Y_{IGBT\text{-}R} + Y_{net} = \begin{bmatrix} Q_1 & Q_2 \\ Q_3 & Q_4 \end{bmatrix}, \tag{12}$$

where $Q_1$ is the ac-port block; $Q_4$ is the dc-port block; $Q_2$ and $Q_3$ represent the coupling between the ac and dc ports. Using the Schur complement, the dc-port reduced characteristic equation is obtained as:

$$Q_4 - Q_3 Q_1^{-1} Q_2 = 0. \tag{13}$$

The equivalent dc-port admittance of the IGBT-R side is therefore defined as:

$$Y_{dceq} = Y_{dc,dc} - Q_3 Q_1^{-1} Q_2. \tag{14}$$

Accordingly, the unit-level dc-port characteristic equation becomes:

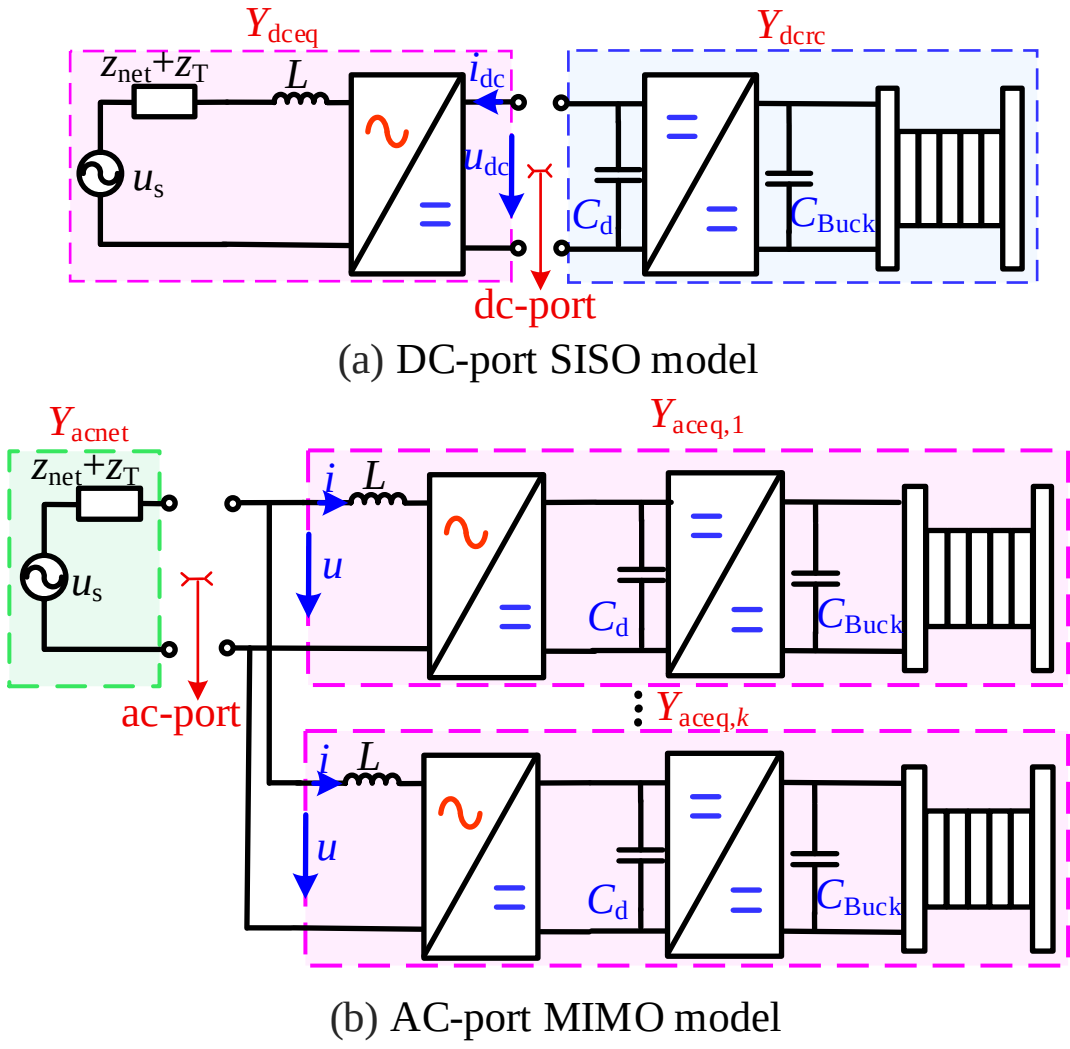

(a) DC-port SISO model

(b) AC-port MIMO model

Fig. 4. Port-reduced equivalent models for unit- and plant-level stability analysis

$$\det(Y_{dceq} + Y_{dcrc}) = 0. \tag{15}$$

As illustrated in Fig. 4(a), $Y_{dceq}$ represents the equivalent admittance of the IGBT-R side reflected to the dc port, whereas $Y_{dcrc}$ represents the dc-link capacitor and the Buck-interfaced AWE stack. The critical frequency $f_c$ is determined by the magnitude intersection

$$\left|Y_{dceq}(j2\pi f_c)\right| = \left|Y_{dcrc}(j2\pi f_c)\right|. \tag{16}$$

The phase difference at the critical frequency is defined as:

$$\Delta\varphi = \left|\angle Y_{dceq}(j2\pi f_c) - \angle Y_{dcrc}(j2\pi f_c)\right|. \tag{17}$$

The unit is regarded as stable when $\Delta\varphi<180°$, critical when $\Delta\varphi=180°$, and unstable when $\Delta\varphi>180°$. This criterion is used in Section III to characterize the dc-port stability boundary under different operating and interface parameters.

### D. Plant-Level AC-Port MIMO Stability Criterion

For plant-level analysis, the dc-side equivalent admittance $Y_{dcrc,k}$ is connected to the dc port of the $k_{th}$ IGBT-R. Using the same block partition as in (12), elimination of the dc port gives:

$$\det\left(Q_{1,k} - Q_{2,k} Q_{4,k}^{-1} Q_{3,k}\right) = 0. \tag{18}$$

Although $Q_{1,k}$ contains the external ac network, it is not included in the unit-equivalent admittance below. Since $Y_{net}$ has no ac- and dc-coupling terms, $Q_{2,k}$ and $Q_{3,k}$ contain only the IGBT-R coupling dynamics. The block $Q_{4,k}$ contains the dc-port admittance of the IGBT-R and $Y_{dcrc,k}$. Thus, the ac-port equivalent admittance of the complete $k_{th}$ unit is

$$Y_{aceq,k} = Y_{IGBT\text{-}R,k}^{2\times2} - Q_{2k} Q_{4k}^{-1} Q_{3k}. \tag{19}$$

Therefore, $Y_{aceq,k}$ retains the effects of the dc-link capacitor, Buck converter, and AWE stack, but does not include the external ac network. When $N$ units are connected in parallel at the PCC, the plant-level ac-port admittance is obtained as:

$$Y_{plant} = \sum_{k=1}^{N} Y_{aceq,k}. \tag{20}$$

The characteristic equation of the plant connected to the external ac network is:

$$\det(Y_{\text{net}}^{2\times 2} + Y_{\text{plant}}) = 0\,. \tag{21}$$

Defining the plant-level return-ratio matrix as:

$$L_{\text{ac}} = Y_{\text{plant}} Z_{\text{ac}}\,, \qquad Z_{\text{ac}} = \left(Y_{\text{net}}^{2\times 2}\right)^{-1}. \tag{22}$$

The plant-level stability is evaluated using the MIMO Nyquist trajectories of $L_{\text{ac}}$. Since the open-loop system has no right-half-plane poles over the operating range considered, plant stability requires the eigenvalue trajectories of $L_{\text{ac}}$ to have no encirclement of $(-1, j0)$. A smaller distance to this point indicates a lower stability margin.

As illustrated in Fig. 4(b), the external ac network interacts with the aggregated admittance of all online IGBT-AWE units. Since $Y_{\text{aceq},k}$ depends on the operating point of each unit, different unit-commitment and power-allocation schemes may produce different Nyquist trajectories under the same total power command. This criterion forms the basis for the plant-level power allocation schemes in Section IV.

## III. Unit-Level DC-Port Stability Characteristics

Fig. 5 applies the dc-port criterion to unit loading, stack temperature, Buck current-loop bandwidth, and dc-link capacitance, with the remaining parameters fixed at the values in Table I. The first two are operating variables, while the latter two are interface parameters. The baseline interface parameters in Table I are deliberately selected near the dc-port stability boundary to reveal parameter-dependent transitions.

### A. Effects of Operating Conditions

Fig. 5(a) presents the dc-port admittances at different AWE current references. At the rated current of 8333 A, $Y_{\text{dceq}}$ and $Y_{\text{dcrc}}$ intersect at 122.3 Hz, where the phase difference is 183.1°. This value exceeds the 180° boundary and indicates an unstable dc-port interaction. Reducing the current reference shifts the magnitude intersection and decreases the phase difference. Lower unit loading therefore increases the dc-port stability margin.

Fig. 5(b) shows the effect of stack temperature at the same current reference, with the steady stack voltage recalculated for each temperature. As the temperature increases from 343 to 373 K, the intersection frequency rises from 119.9 to 125.6 Hz, while the phase difference increases from 179.5° to 188.3°. The operating point crosses the stability boundary between 343 and 353 K. Increasing the stack temperature therefore reduces the dc-port stability margin. This comparison isolates the electrical stability effect, while the stack voltage limit is included in the operating region constructed in Section IV.

### B. Effects of Interface Parameters

Fig. 5(c) examines the Buck current-loop bandwidth $f_{\text{bcl}}$. Changing the bandwidth mainly reshapes $Y_{\text{dcrc}}$, while the rectifier operating point and $Y_{\text{dceq}}$ remain nearly unchanged. Within the tested range, increasing the bandwidth reduces the phase difference at the magnitude intersection and moves the unit away from the critical condition. The Buck bandwidth should therefore be selected by considering both current response and dc-port stability.

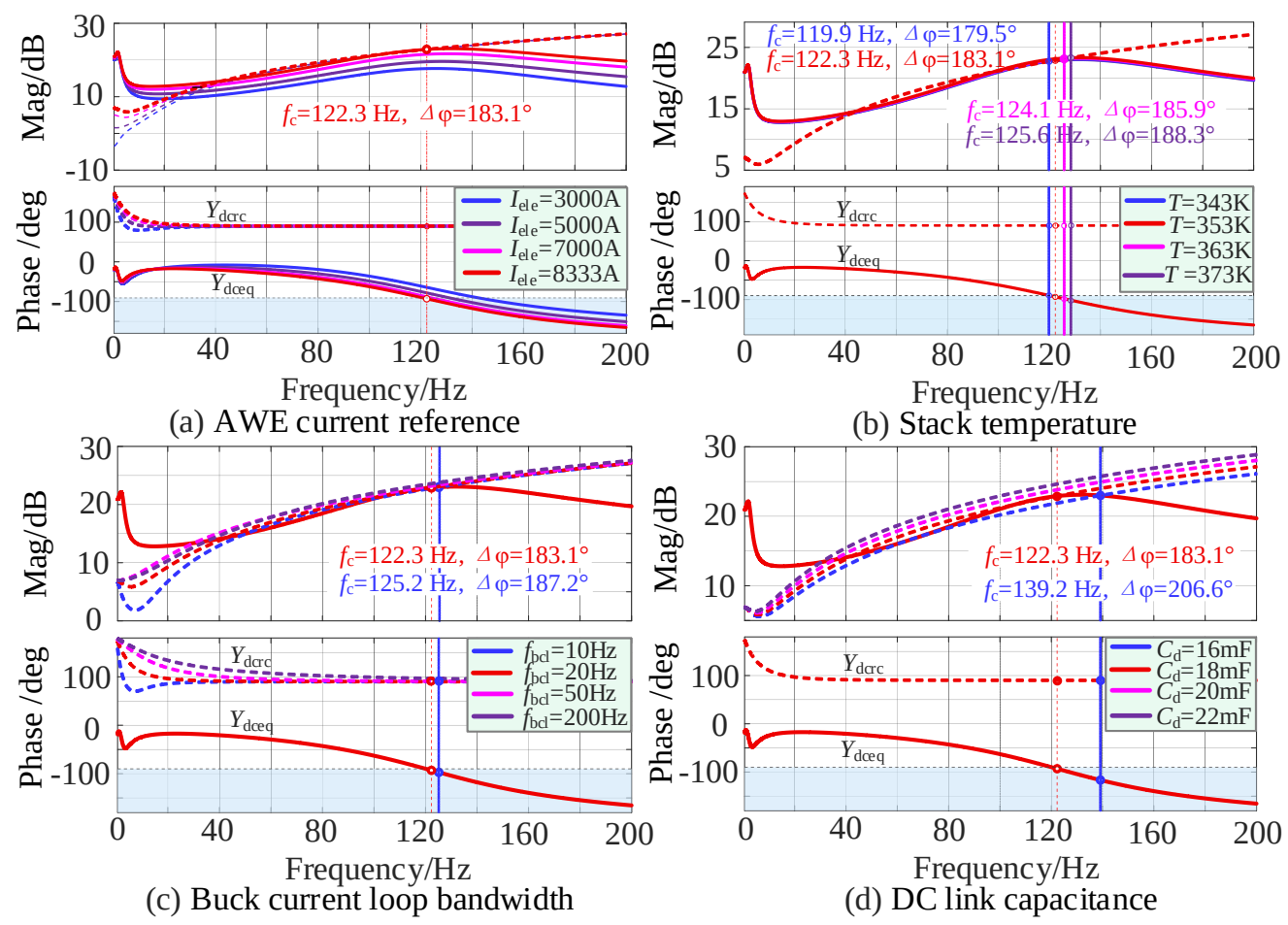

Fig. 5. Unit-level dc-port stability characteristics

Fig. 5(d) presents the effect of the dc-link capacitance. Increasing $C_{\text{d}}$ raises the magnitude of $Y_{\text{dcrc}}$ near the critical frequency and shifts the intersection toward a lower frequency. For example, increasing $C_{\text{d}}$ from 16 to 18 mF reduces the phase difference from 206.6° to 183.1°. A larger dc-link capacitance therefore improves the stability margin within the considered parameter range.

## IV. Stability-Oriented Operating Regions and Power Allocation Dispatch

### A. Stability-Oriented Dispatch Framework

Large ReP2H plants use EMS to coordinate unit commitment and power allocation subject to loading, temperature, and gas impurity constraints under fluctuating renewable power supply [16], [17], [18]. These constraints ensure production feasibility but do not account for operating-point-dependent electrical stability. The proposed framework therefore uses offline operating and design regions to screen unit loading, temperature, Buck bandwidth, and dc-link capacitance, followed by plant-level MIMO Nyquist assessment of candidate unit-commitment and power-allocation schemes. Schedules with insufficient margins are rejected or adjusted before dispatch, and the resulting lookup maps avoid repeated online frequency-domain analysis.

### B. Unit-Level DC-Port Operating and Design Regions

Based on the dc-port stability criterion in Section II-C, the phase difference between $Y_{\text{dceq}}$ and $Y_{\text{dcrc}}$ at their magnitude intersection is used to characterize the proximity of an operating point to the dc-port oscillatory boundary. The corresponding dc-port phase margin is defined as

$$PM_{\text{dc}} = 180^{\circ} - \Delta\varphi\,, \tag{23}$$

where $\Delta\varphi$ is the phase difference at the magnitude intersection. An operating point without a magnitude intersection is regarded as stable. Otherwise, the operating points are classified as:

$$\begin{cases} \text{Stable, } PM_{\text{dc}} > 20^{\circ} & \text{Low risk, } 10^{\circ} < PM_{\text{dc}} \le 20^{\circ} \\ \text{High risk, } 0^{\circ} < PM_{\text{dc}} \le 10^{\circ} & \text{Unstable, } PM_{\text{dc}} \le 0^{\circ} \end{cases} \tag{24}$$

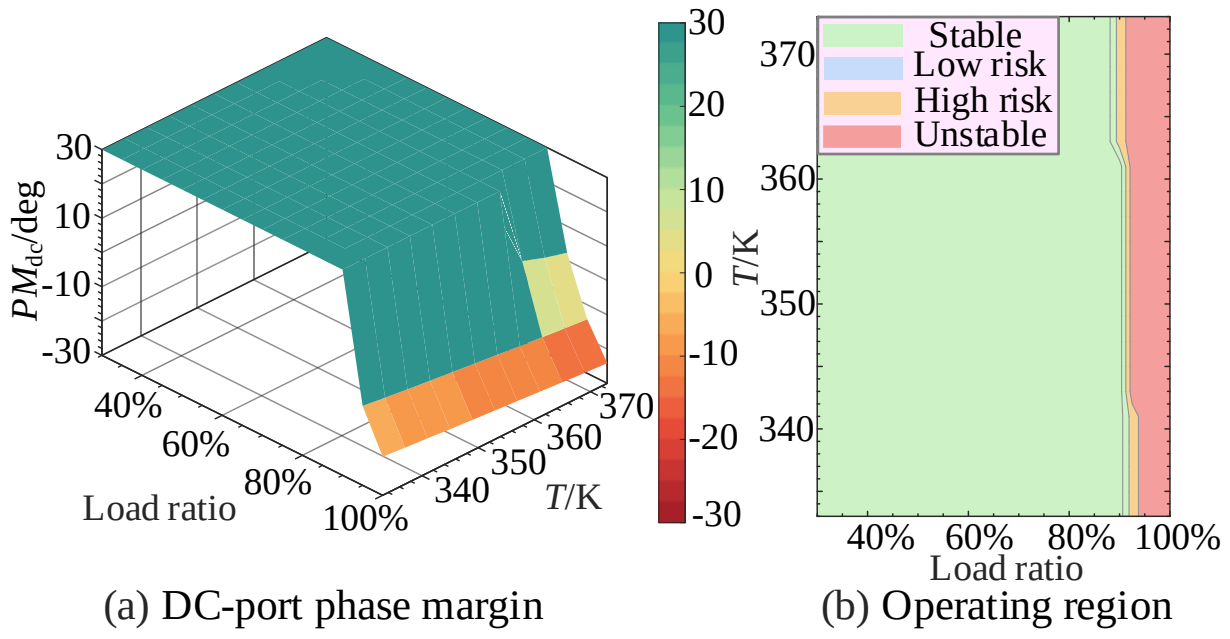

(a) DC-port phase margin (b) Operating region

Fig. 6. Unit-level dc-port operating region under fixed interface parameters

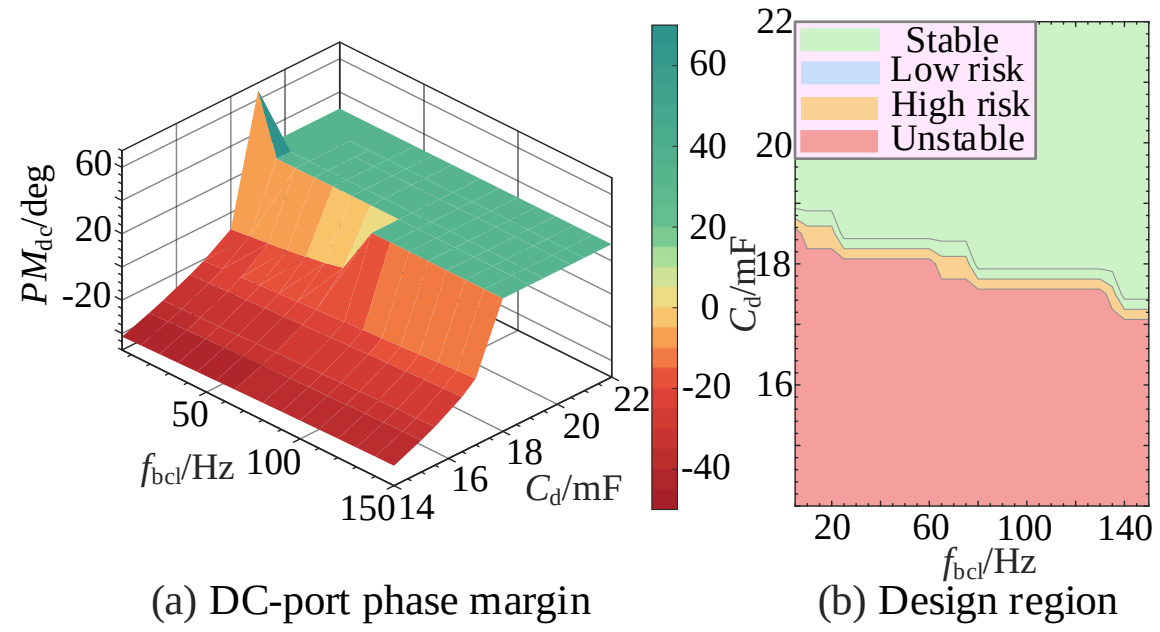

(a) DC-port phase margin (b) Design region

Fig. 7. Unit-level dc-port design region at the rated operating condition

The four levels classify operating points according to their proximity to the dc-port oscillatory boundary. For this case study, the 10° and 20° thresholds serve as adjustable engineering margins rather than universal limits. Low-risk and high-risk points remain stable but are increasingly sensitive to parameter variations and operating point changes, while points with $PM_{dc} \leq 0°$ are classified as unstable.

Fig. 6 presents the unit-level operating region in the load ratio and stack temperature plane, with the Buck current-loop bandwidth and dc-link capacitance fixed at 50 Hz and 18 mF, respectively. As shown in Fig. 6(a), the dc-port phase margin remains high over most low- and medium-load conditions, while it decreases rapidly as the load approaches its rated value. Increasing the stack temperature further reduces the margin near the high-load boundary.

The categorized map in Fig. 6(b) shows that most operating points are stable, whereas the high-risk and unstable regions are concentrated near the rated load. For visualization, points without a magnitude intersection are displayed at $PM_{dc}$=30°. The resulting map allows the EMS to determine the admissible unit loading according to the measured stack temperature.

Fig. 7 evaluates the effects of the $f_{bcl}$ and $C_d$ at the rated current and 353 K. Fig. 7(a) shows that the dc-link capacitance has a dominant influence on the phase margin. Increasing $C_d$ increases the dc-port stability margin, while increasing the Buck bandwidth reduces the minimum capacitance required for an adequate margin within the investigated range.

The categorized design map in Fig. 7(b) identifies the combinations that satisfy the prescribed margin levels. Since the Buck bandwidth and dc-link capacitance are control and circuit parameters, respectively, Fig. 7 is used for controller and interface design rather than online dispatch.

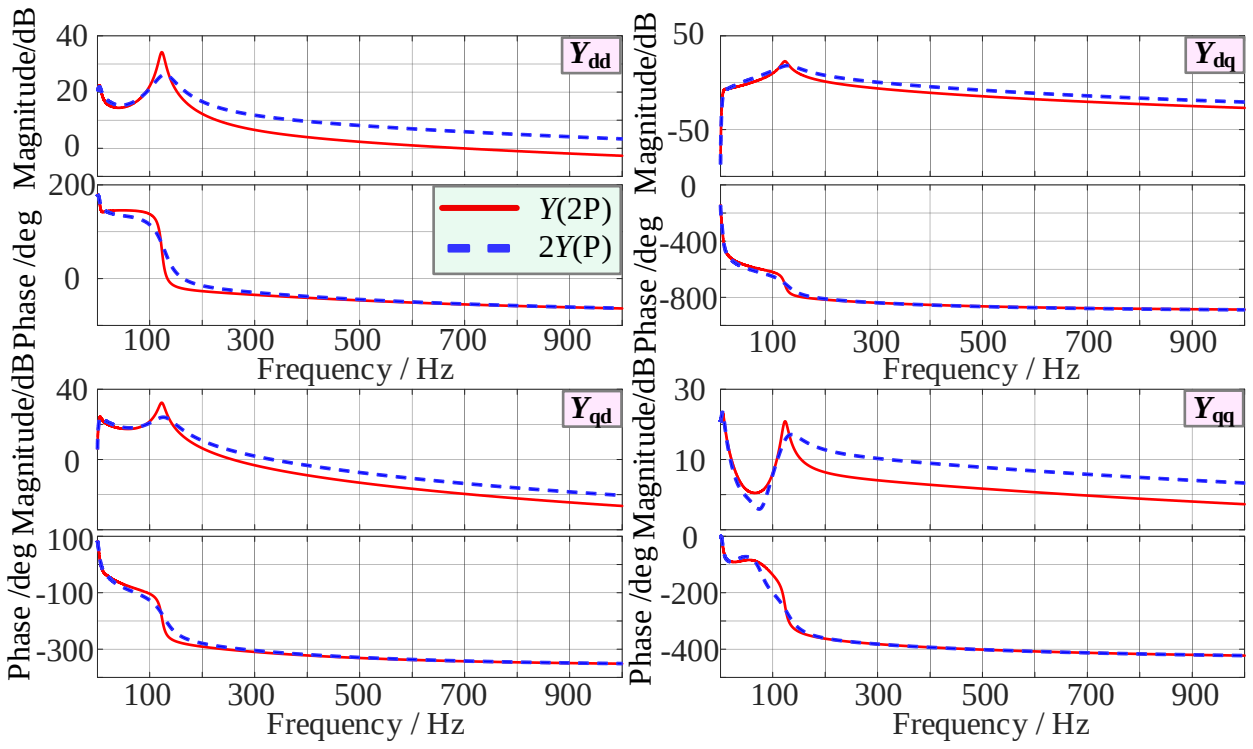

Fig. 8. Comparison between $Y_{aceq}(2P)$ and $2Y_{aceq}(P)$ under the same total power command

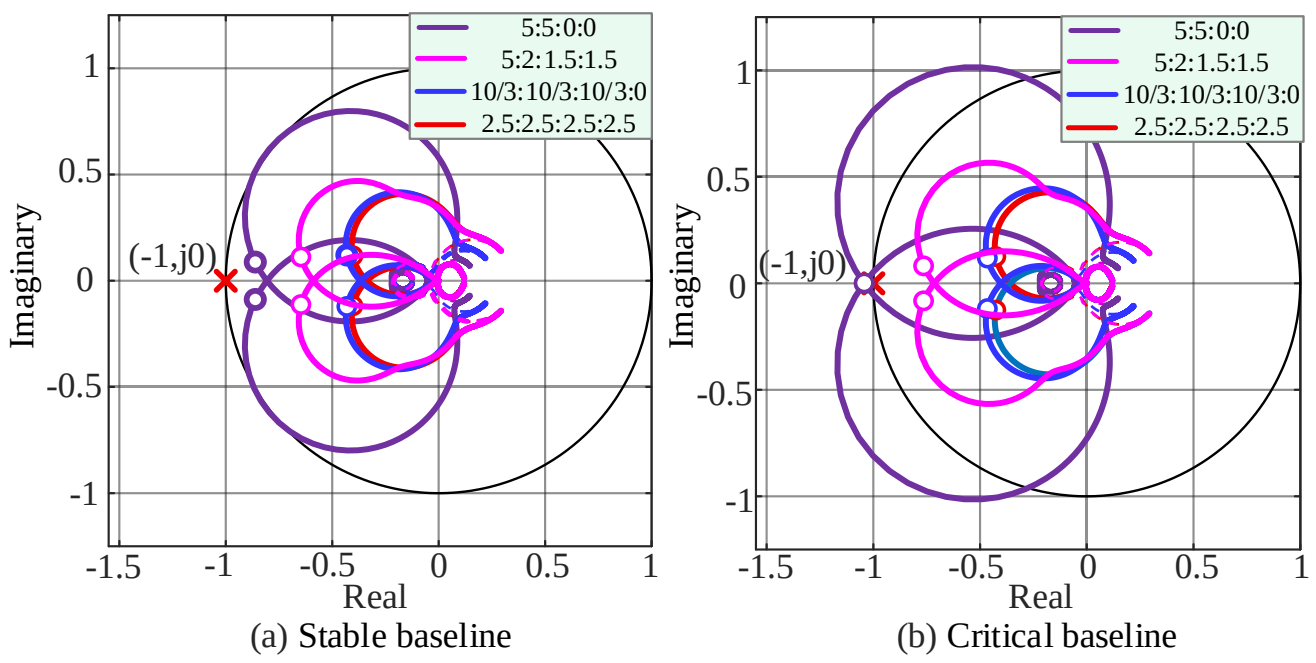

Fig. 9. MIMO Nyquist trajectories of a 4-in-1 AWE module under different power allocation schemes with a total power command of 10 MW

The two maps provide complementary information. Fig. 6 defines the admissible operating range after the interface parameters have been selected, whereas Fig. 7 supports the offline coordination of the Buck controller and dc-link capacitor. Both maps can be generated through offline scanning and stored as lookup tables for Buck-controller/dc-link capacitor coordination and multi-ELZ dispatch.

### *C. Plant-Level Power Allocation Feasible Region*

At the plant level, the same hydrogen-production command can be achieved through different unit-commitment and power allocation schemes. Since the ac-port equivalent admittance of each IGBT-AWE unit depends on its own operating point, the aggregated plant admittance is not proportional to the total power command. Fig. 8 compares the ac-port admittance $Y_{aceq}(2P)$ of one unit operating at power $2P$ with the aggregated admittance $2Y_{aceq}(P)$ of two identical units operating at power $P$. The observed difference confirms that:

$$Y_{aceq}(2P) \neq 2Y_{aceq}(P) \tag{25}$$

and therefore one heavily loaded unit is not electrically equivalent to two lightly loaded units under the same total power command.

A representative 4-in-1 ELZ module comprising four independently powered stacks is considered. In this configuration, four stacks share the gas separation and lye circulation equipment within one balance of plant (BoP) system, while each stack is supplied by an independent rectifier and can be controlled separately. Each IGBT-AWE unit is rated at 5 MW, and the total power command is fixed at 10 MW.

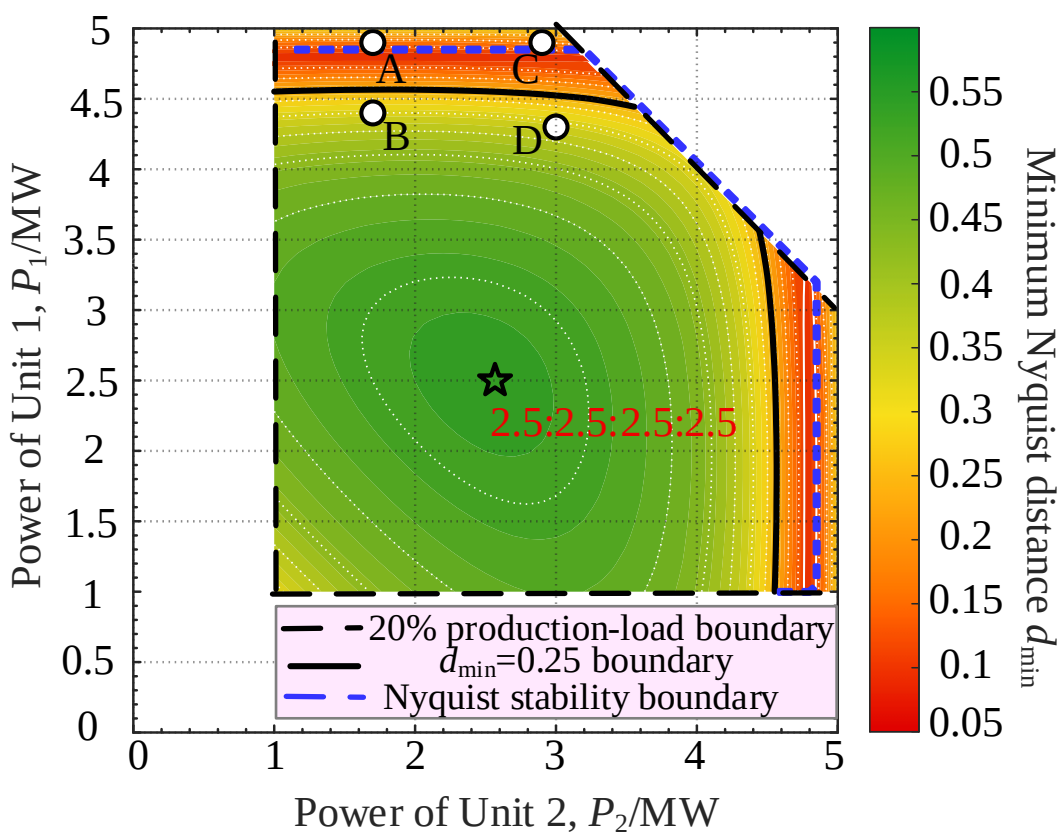

Fig. 10. Representative plant-level dispatch map under a fixed 10 MW command

Fig. 9 compares the MIMO Nyquist eigenvalue trajectories of representative power allocations under the same 10 MW command. For both the stable and critical baselines, the concentrated allocation 5:5:0:0 MW produces the trajectory closest to the critical point (-1, $j$0). Distributing the command among more units moves the trajectories away from the critical point, while equal allocation provides the largest margin among the examined cases. These results show that satisfying the total power command and unit capacity limits alone does not guarantee adequate plant stability.

Based on the unit-level operating regions and the plant-level MIMO criterion, the set of acceptable schedules is defined as:

$$\mathcal{F}_{\text{sch}} = \left\{ (b,P) \;\middle|\; \begin{array}{c} b_k \in \{0,1\}, \qquad k=1,\dots,N \\ \sum_{k=1}^{N} P_k = P_{\text{cmd}}, b_k P_k^{\min} \le P_k \le b_k P_k^{\max} \\ P_k^{\min} = 0.2 P_k^{\text{rated}}, k=1,\dots,N, \\ (P_k, T_k) \in \mathcal{F}_{\text{dc},k}, \qquad \forall k \text{ with } b_k = 1, \\ N_{\text{enc}}\left[L_{\text{ac}}(b,P)\right] = 0 \\ d_{\min}(b,P) \ge d_{\text{req}} \end{array} \right\}, \quad (26)$$

where $\mathcal{F}_{\text{sch}}$ denotes the set of acceptable schedules, and $N$ is the number of IGBT-AWE units. The vector $b$=[$b_1$,...,$b_N$]$^T$ contains the unit states. The binary variable $b_k$=1 denotes online production, while $b_k$=0 denotes the zero-power state. The vector $P$=[$P_1$,...,$P_N$]$^T$ contains unit power commands, and $P_{\text{cmd}}$ is the total production command. The parameters $P_k^{\min}$ and $P_k^{\max}$ denote the power limits of the $k_{\text{th}}$ online unit. To account for the minimum load imposed by the hydrogen-in-oxygen (HTO) impurity limit, the minimum production power is set to 20% of the rated unit power. The stack temperature is denoted by $T_k$, and $\mathcal{F}_{\text{dc},k}$ is its unit-level dc-port operating region. The function $N_{\text{enc}}$[.] gives the number of encirclements of the critical point (-1, $j$0) by the MIMO Nyquist eigenvalue trajectories. The quantities $d_{\min}$ and $d_{\text{req}}$ denote the minimum Nyquist distance and its prescribed threshold, respectively.

The plant return ratio depends on the operating state and power command of each unit and is expressed as:

$$L_{\text{ac}}(b,P,s) = \left[ \sum_{k=1}^{N} \left( b_k Y_{\text{aceq},k}^{\text{on}}(P_k,s) + (1-b_k) Y_{\text{aceq},k}^{\text{st}}(s) \right) \right] Z_{\text{ac}}(s), \quad (27)$$

where $Y_{\text{aceq},k}^{\text{on}}$ and $Y_{\text{aceq},k}^{\text{st}}$ are the equivalent ac port admittances of the $k_{\text{th}}$ unit at the production and zero-power operating points, respectively. The matrix $Z_{\text{ac}}$ is the external ac network impedance defined in (22).

The minimum distance between the MIMO Nyquist eigenvalue trajectories and the critical point is:

$$d_{\min}(b,P) = \min_{\omega,m} \left| 1 + \lambda_m \left[ L_{\text{ac}}(b,P,j\omega) \right] \right|, \quad (28)$$

where $\lambda_m$[.] denotes the $m_{\text{th}}$ eigenvalue of the return ratio matrix; and $\omega$ is the angular frequency. A schedule is accepted when all production units remain within their dc-port operating regions, the Nyquist trajectories do not encircle the critical point, and $d_{\min}$ remains above $d_{\text{req}}$.

Fig. 10 extends the discrete comparisons in Fig. 9 to a representative continuous dispatch slice under the fixed 10 MW command. We vary the powers of Units 1 and 2 independently and assign the remaining power equally to Units 3 and 4:

$$P_3 = P_4 = \frac{10 - P_1 - P_2}{2}, \; P_k \in \{0\} \cup [1,5] \text{ MW}, \; k=1,2,3,4 \,. \quad (29)$$

In Fig. 10, $P_k$=0 represents the zero-power state, while 1≤$P_k$≤5 MW denotes online production. Operating points with $0 < P_k < 1$ MW are excluded by the 20% minimum load constraint. The color represents $d_{\min}$, and the black contour marks $d_{\min}$=0.25. The blue dashed line denotes the Nyquist stability boundary determined from the encirclement condition in (26). Equal sharing provides the largest margin, while concentrated allocations reduce $d_{\min}$, consistent with Fig. 9. To further examine the predicted stability boundary, four representative operating points are selected. Their power allocations ([$P_1$, $P_2$, $P_3$, $P_4$]) are A= [4.90, 1.70, 1.70, 1.70] MW, B= [4.40, 1.70, 1.95, 1.95] MW, C= [4.90, 2.90, 1.10, 1.10] MW, and D= [4.30, 3.00, 1.35, 1.35] MW, respectively. Cases A and C are predicted to be unstable, whereas Cases B and D are predicted to be stable. These four cases are subsequently examined by plant-level HIL tests in Section V. The map therefore identifies schedules that satisfy the production command and unit power limits but provide insufficient electrical stability margin.

For larger plants, the unit-level operating regions can first eliminate infeasible operating points. Plant-level MIMO analysis can then evaluate the remaining unit commitment and power sharing schemes, reducing the required number of frequency-domain evaluations.

## V. Hardware-in-the-Loop Validation

### A. HIL Platform and Test Setup

Fig. 11 shows the HIL platform used to validate the proposed stability-oriented operating regions and power allocation schemes. The power stage and the interconnected plant model are configured in StarSim and executed on an MT8020 real-time simulator with a time step of 1 µs. The simulator reproduces the electrical dynamics of the IGBT-AWE unit, including the ac-grid interface, dc-link capacitor, Buck converter, and AWE stack. Voltage and current signals are exchanged with the DSP and FPGA controller, which returns the corresponding pulse-width modulation (PWM) commands to the real-time model. This closed-loop implementation preserves the converter control and dc-link interactions required by the stability analysis.

### B. Unit-Level Validation of the DC-Port Stability Boundary

Fig. 12 confirms the normal steady-state operation of the IGBT-AWE unit, with sinusoidal ac voltage and current, regulated dc-link voltage, and stable AWE current and stack voltage.

The operating point is then shifted toward the dc-port stability boundary. As shown in Fig. 13(a), sustained oscillations appear simultaneously in the dc-link voltage, AWE current, and stack voltage. The enlarged waveforms in Fig. 13(b) exhibit the same dominant oscillatory component, and the fast Fourier transform (FFT) identifies a frequency of approximately 121 Hz. This agrees well with the critical frequency predicted by the dc-port admittance analysis, validating the magnitude-intersection and phase-difference criterion in Section II-C.

Fig. 14 further verifies the parameter trends obtained in Section III. Increasing the Buck current-loop bandwidth from 20 to 100 Hz suppresses the oscillation, consistent with the predicted increase in dc-port stability margin. In contrast, increasing the stack temperature from 353 to 373 K enlarges the oscillation, confirming the adverse effect of temperature on the stability boundary. The HIL results therefore reproduce the effects of both interface and operating parameters identified by the frequency-domain analysis.

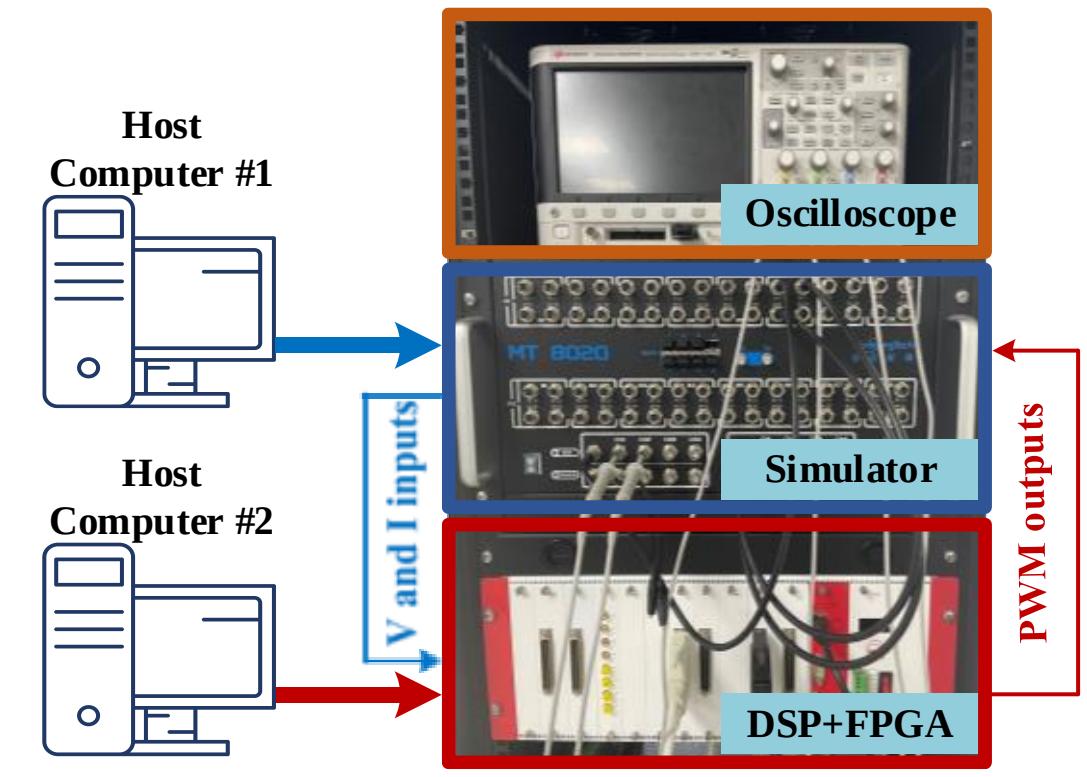

Fig. 11. HIL platform for validating the IGBT–AWE electrolysis unit

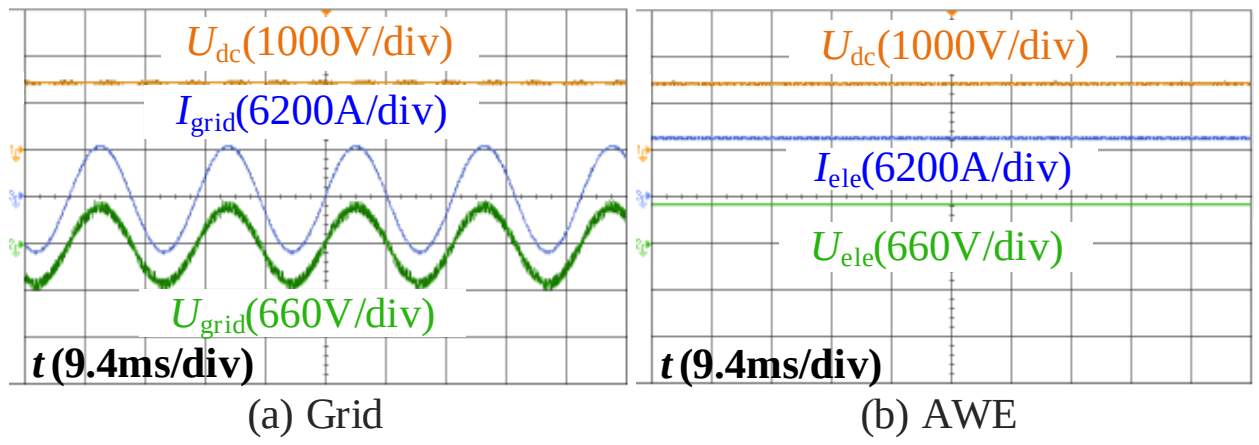

(a) Grid (b) AWE

Fig. 12. Steady-state HIL waveforms of one IGBT-AWE unit

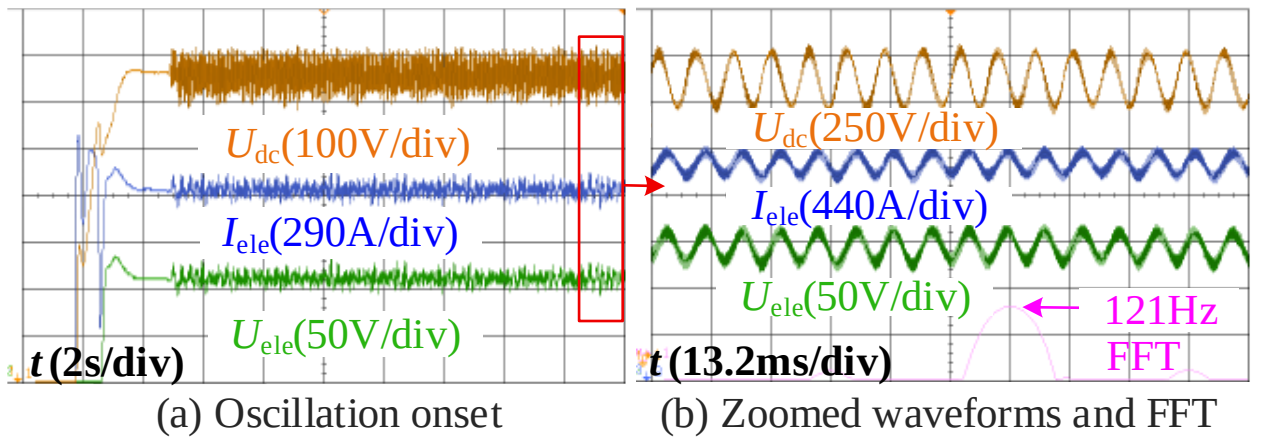

(a) Oscillation onset (b) Zoomed waveforms and FFT

Fig. 13. HIL validation of the dc-port oscillation near the stability boundary

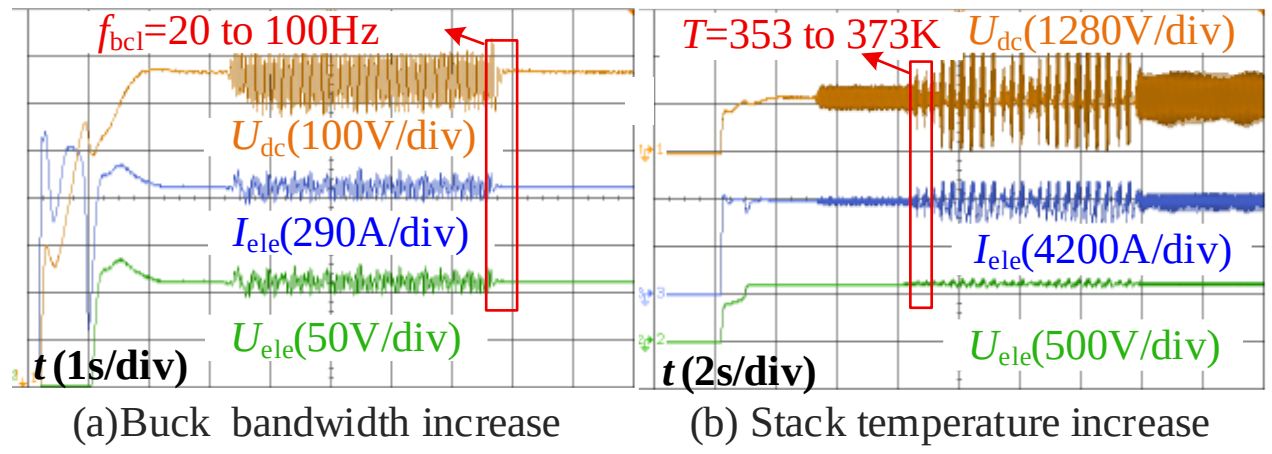

(a)Buck bandwidth increase (b) Stack temperature increase

Fig. 14. HIL responses under parameter variations

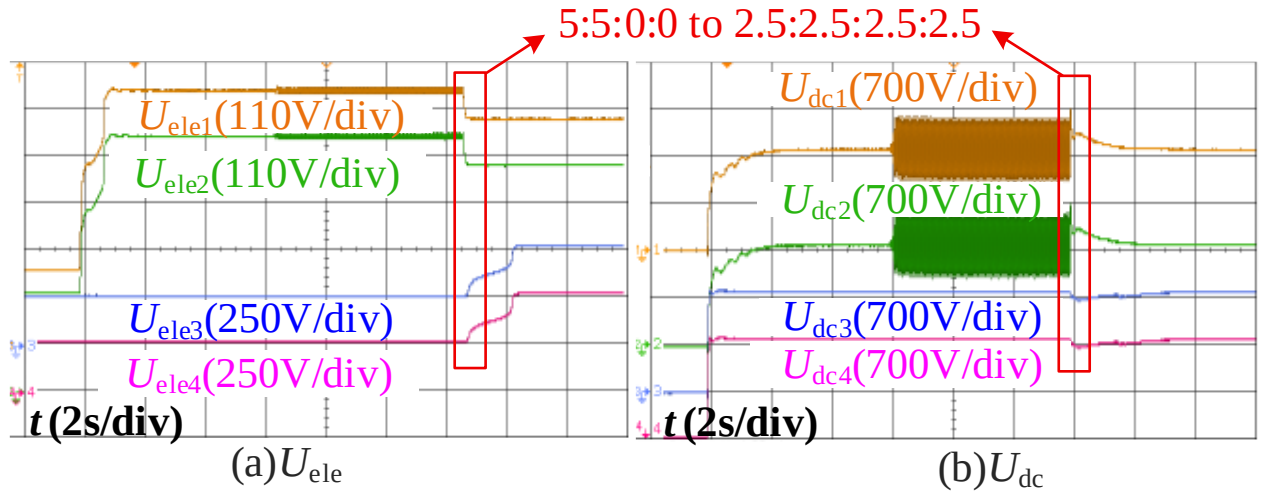

(a)$U_{ele}$ (b)$U_{dc}$

Fig. 15. Plant-level HIL response during power redistribution from 5:5:0:0 to 2.5:2.5:2.5:2.5 MW

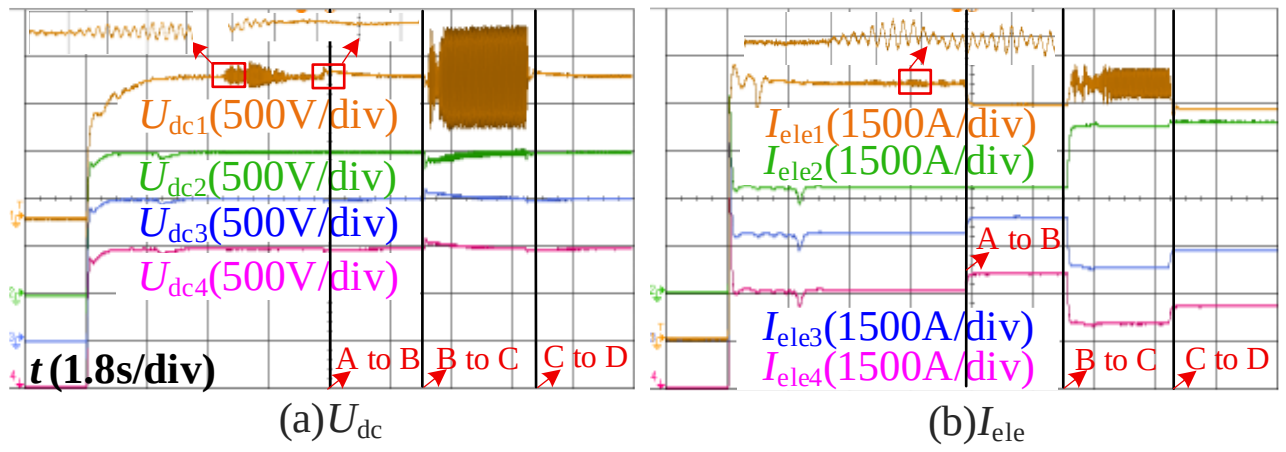

(a)$U_{dc}$ (b)$I_{ele}$

Fig. 16. Plant-level HIL validation of representative Cases A to D across the predicted stability boundary

### C. Validation of Plant-Level Power Allocation

Fig. 15 validates the stability-oriented power allocation rule derived from Figs. 9 and 10. Under the initial 5:5:0:0 MW allocation, Units 1 and 2 operate at rated power beyond the unit-level dc-port stability boundary identified in Section III, while Units 3 and 4 remain at zero power. Sustained oscillations appear in the dc-link voltages and stack voltages of the loaded units, consistent with the unit-level stability prediction.

The power command is then redistributed to 2.5:2.5:2.5:2.5 MW. All four units consequently operate at half load, which lies within the stable unit operating region, while the plant-level MIMO analysis also predicts an increased stability margin. The measured dc-link oscillations decay after redistribution and the four units converge to stable operating points. The HIL results therefore confirm that power redistribution can increase the stability margin and suppress oscillations under the same total production command.

Fig. 16 further evaluates the four nonuniform schedules marked in Fig. 10. Cases A and C are predicted to be unstable, whereas Cases B and D are stable. The enlarged segments in Fig. 16 highlight the relatively weak oscillatory component in Case A, while a pronounced sustained oscillation is observed in Case C. After the operating point is changed from A to B and from C to D, the oscillations decay and the responses converge to stable operating points. The alternating responses of Cases A to D agree with the stability boundary predicted in Fig. 10 and confirm that power allocation can change the plant-level stability margin under the same total production command.

## VI. Conclusion

This paper develops a stability-oriented framework that links the dc-port dynamics of individual IGBT-AWE units with the aggregated ac-port behavior of multi-ELZ ReP2H plants. The analytical and HIL results yield the following findings.

1) At the unit level, higher stack temperature and heavier loading reduce the dc-port stability margin. Increasing the dc-link capacitance and Buck current-loop bandwidth improves stability. The resulting design and operating regions identify parameter combinations and operating points close to the oscillation boundary.
2) At the plant level, the aggregated ac-port admittance depends on unit commitment and power allocation. The same total production power can therefore produce different stability margins. Balanced power allocation among more online units generally increases the stability margin, whereas concentrated allocation moves the operating point closer to the plant-level stability boundary.
3) The proposed operating and design regions translate the identified stability margins into practical constraints for interface design and production dispatch. HIL tests verify the predicted stability boundaries and confirm that redistributing power from low-margin or unstable allocations to higher-margin allocations can suppress oscillations without changing the total production command.

Future work will integrate the proposed stability-oriented operation rule into the plant EMS to coordinately optimize production cost, startup and standby losses, degradation, and process constraints. Large plants with heterogeneous ELZs and rectifiers will be addressed through hierarchical screening and reduced-order aggregation.

## Appendix A

Grid voltage perturbations produce the angle deviation $\Delta\theta$, and linearization gives the PLL coupling matrices:

$$\begin{bmatrix}\Delta u_d^c\\ \Delta u_q^c\\ \Delta u_{dc}\end{bmatrix}=\underbrace{\begin{bmatrix}1 & G_{PLL}U_q^s & 0\\ 0 & 1-G_{PLL}U_d^s & 0\\ 0 & 0 & 1\end{bmatrix}}_{G^u{}_{PLL}}\begin{bmatrix}\Delta u_d^s\\ \Delta u_q^s\\ \Delta u_{dc}\end{bmatrix},$$

$$\begin{bmatrix}\Delta i_d^c\\ \Delta i_q^c\\ \Delta i_{dc}\end{bmatrix}=\underbrace{\begin{bmatrix}0 & G_{PLL}I_q^s & 0\\ 0 & -G_{PLL}I_d^s & 0\\ 0 & 0 & 0\end{bmatrix}}_{G^i{}_{PLL}}\begin{bmatrix}\Delta u_d^s\\ \Delta u_q^s\\ \Delta u_{dc}\end{bmatrix}, \quad \text{(A1)}$$

$$\begin{bmatrix}\Delta d_d^c\\ \Delta d_q^c\\ \Delta d_{dc}\end{bmatrix}=\underbrace{\begin{bmatrix}0 & -G_{PLL}D_q^s & 0\\ 0 & G_{PLL}D_d^s & 0\\ 0 & 0 & 0\end{bmatrix}}_{G^d{}_{PLL}}\begin{bmatrix}\Delta u_d^s\\ \Delta u_q^s\\ \Delta u_{dc}\end{bmatrix}, G_{PLL}=\frac{G_{PI_PLL}}{s+U_d^sG_{PI_PLL}}$$

where the superscript $c$ denotes the controller $dq$-frame; $G_{PI_PLL} = K_{P_PLL}+ K_{i_PLL}/s$.

The dc-link voltage loop gives the $d$-axis current reference:

$$\begin{bmatrix}\Delta i_{dref}\\ 0\\ 0\end{bmatrix}=-\underbrace{\begin{bmatrix}0 & 0 & G_{PI_DVC}\\ 0 & 0 & 0\\ 0 & 0 & 0\end{bmatrix}}_{H_u}\begin{bmatrix}\Delta u_d^s\\ \Delta u_q^s\\ \Delta u_{dc}\end{bmatrix}, \quad \text{(A2)}$$

where $G_{PI_DVC} = K_{P_DVC}+ K_{i_DVC}/s$.

Linearizing the power balance and modulation relation yields:

$$\begin{bmatrix}\Delta u_{abd}^s\\ \Delta u_{abq}^s\\ \Delta i_{dc}\end{bmatrix}=\underbrace{\begin{bmatrix}0.5U_{dc} & 0 & 0.5D_d\\ 0 & 0.5U_{dc} & 0.5D_q\\ -0.75I_d & -0.75I_q & 0\end{bmatrix}}_{G_{m1}}\begin{bmatrix}\Delta d_d^s\\ \Delta d_q^s\\ \Delta u_{dc}\end{bmatrix} + \underbrace{\begin{bmatrix}0 & 0 & 0\\ 0 & 0 & 0\\ -0.75D_d & -0.75D_q & 0\end{bmatrix}}_{G_{m2}}\begin{bmatrix}\Delta i_d^s\\ \Delta i_q^s\\ \Delta i_{dc}\end{bmatrix}, \quad \text{(A3)}$$

where the first two variables on the left of (A3) are the $d$ and $q$ components of the rectifier bridge output voltage. Capital letters denote steady-state quantities. In particular, $D_d$ and $D_q$ are the steady-state $d$ and $q$ modulation ratios.

The ac current loop, including the $dq$ cross-coupling terms, is written as:

$$\begin{bmatrix}\Delta u_{abdref}^c\\ \Delta u_{abqref}^c\\ \Delta u_{dc}\end{bmatrix}=\begin{bmatrix}\Delta u_d^c\\ \Delta u_q^c\\ \Delta u_{dc}\end{bmatrix}-\underbrace{\begin{bmatrix}G_{PI_ACC} & 0 & 0\\ 0 & G_{PI_ACC} & 0\\ 0 & 0 & 0\end{bmatrix}}_{H_i}\begin{bmatrix}\Delta i_{dref}-\Delta i_d^c\\ \Delta i_{qref}-\Delta i_q^c\\ 0-\Delta i_{dc}\end{bmatrix} + \underbrace{\begin{bmatrix}0 & \omega_0 L & 0\\ -\omega_0 L & 0 & 0\\ 0 & 0 & 0\end{bmatrix}}_{G_{ol}}\begin{bmatrix}\Delta i_d^c\\ \Delta i_q^c\\ \Delta i_{dc}\end{bmatrix}, \quad \text{(A4)}$$

where $G_{PI_ACC} = K_{P_ACC} + K_{i_ACC}/s$, and $\omega_0$ is the fundamental angular frequency.

The modulation reference is normalized by the dc-link voltage as:

$$\begin{bmatrix}\Delta d_{dref}^c\\ \Delta d_{qref}^c\\ \Delta u_{dc}\end{bmatrix}=\underbrace{\begin{bmatrix}0.5U_{dc} & 0 & 0.5D_d\\ 0 & 0.5U_{dc} & 0.5D_q\\ 0 & 0 & 1\end{bmatrix}^{-1}}_{G_{m3}}\begin{bmatrix}\Delta u_{abdref}^c\\ \Delta u_{abqref}^c\\ \Delta u_{dc}\end{bmatrix}, \quad \text{(A5)}$$

With the PWM delay $T_d$, the duty cycle becomes:

$$\begin{bmatrix}\Delta d_d^s\\ \Delta d_q^s\\ \Delta u_{dc}\end{bmatrix}=\underbrace{\begin{bmatrix}e^{-T_d s} & 0 & 0\\ 0 & e^{-T_d s} & 0\\ 0 & 0 & 1\end{bmatrix}}_{G_d}\begin{bmatrix}\Delta d_{dref}^c\\ \Delta d_{qref}^c\\ \Delta u_{dc}\end{bmatrix}, \quad \text{(A6)}$$

The ac filter gives the port voltage and current relation:

$$\begin{cases}\underbrace{\begin{bmatrix}F_1 & F_2 & 0\\ -F_2 & F_1 & 0\\ 0 & 0 & 0\end{bmatrix}}_{G_{rl1}}\begin{bmatrix}\Delta u_d^s\\ \Delta u_q^s\\ \Delta u_{dc}\end{bmatrix}-\underbrace{\begin{bmatrix}F_1 & F_2 & 0\\ -F_2 & F_1 & 0\\ 0 & 0 & -1\end{bmatrix}}_{G_{rl2}}\begin{bmatrix}\Delta u_{abd}^s\\ \Delta u_{abq}^s\\ \Delta i_{dc}\end{bmatrix}=\begin{bmatrix}\Delta i_d^s\\ \Delta i_q^s\\ \Delta i_{dc}\end{bmatrix},\\ F_1=\dfrac{R_f+sL_f}{(R_f+sL_f)^2+\omega_0^2L_f^2}, F_2=\dfrac{\omega_0 L_f}{(R_f+sL_f)^2+\omega_0^2L_f^2}\end{cases} \quad \text{(A7)}$$

where $R_f$ and $L_f$ are the filter resistance and inductance, respectively.

Combining (A1) to (A7) with the rectifier power balance equations gives the three-port admittance matrix:

$$\underbrace{(G_{\mathrm{rl2}}{}^{-1} + G_{\mathrm{Nt}})^{-1}(G_{\mathrm{rl2}}{}^{-1}G_{\mathrm{rl1}} - G_{\mathrm{Mt}})}_{Y_{\mathrm{IGBT\text{-}R}}}\Delta u^{\mathrm{s}} = \Delta i^{\mathrm{s}}. \tag{A8}$$

where

$$\begin{cases} G_{\mathrm{Mt}} = G_{\mathrm{d}}G_{\mathrm{m1}}(G_{\mathrm{m3}}(H_{\mathrm{i}}H_{\mathrm{u}} + G_{\mathrm{PLL}}^{u} + G_{\mathrm{PLL}}^{i}(H_{\mathrm{i}} + G_{\omega\mathrm{l}})) + G_{\mathrm{PLL}}^{d}) \\ G_{\mathrm{Nt}} = G_{\mathrm{d}}G_{\mathrm{m1}}G_{\mathrm{m3}}(H_{\mathrm{i}} + G_{\omega\mathrm{l}}) + G_{\mathrm{m2}} \end{cases} \tag{A9}$$

The matrices in (A9) collect the voltage-driven and current-driven dynamics of the IGBT-R. Equation (A8) represents the IGBT-R three-port admittance. The dc-link capacitor and the Buck-interfaced AWE stack are connected externally to the dc port through the dc-side admittance defined in (3).